\documentclass{article}
\usepackage{latexsym}
\usepackage{graphics}
\usepackage{amsfonts}
\usepackage{graphicx}
\usepackage{epstopdf}
\newtheorem{thm}{Theorem}

\newtheorem{lem}{Lemma}
\newtheorem{defn}{Definition}
\newtheorem{cor}{Corollary}

\begin{document}

\title{Weight systems for Milnor invariants}
\author{Blake Mellor\\
			Mathematics Department\\
			Loyola Marymount University\\
			Los Angeles, CA  90045-2659\\
   {\it  bmellor@lmu.edu}}
\date{}
\maketitle

\begin{abstract}
We use Polyak's skein relation to give a new proof that Milnor's string link invariants $\mu_{12...n}$ are
finite type invariants, and to develop a recursive relation for their associated weight
systems.  We show that the obstruction to the triviality of these weight systems is the presence of a
certain kind of spanning tree in the intersection graph of a chord diagram.
\end{abstract}

\section{Introduction} \label{S:intro}

Milnor's $\mu$-invariants \cite{mi} are important invariants of link homotopy.  Unfortunately, they are only
well-defined modulo a complicated indeterminacy, giving the more common $\bar{\mu}$-invariants.  Habegger and
Lin \cite{hl} noticed that the indeterminacy disappears if the $\mu$-invariants are viewed as invariants of
{\it string links} rather than links, and proved that the $\mu$-invariants classify string links up to
homotopy.  These invariants are also part of the theory of finite-type invariants: it is known that the
invariant $\mu_{i_1i_2...i_r,j}$ is an invariant of type $r$ for string links \cite{bn2, li}.  

Every finite type invariant gives rise to a weight system on chord diagrams \cite{bn}, and these weight systems
can be very useful in understanding the associated invariants, as in Bar-Natan and Garoufalides' proof
of the Melvin-Morton-Rozansky Conjecture \cite{bg}.  To date, the weight systems for the Milnor invariants have
only been described in terms of unitrivalent graphs \cite{bn2, hm}, rather than directly in terms of chord
diagrams.  Recently, Polyak \cite{po} has proven a skein relation for the $\mu$-invariants, analogous to the
well-known skein relations for the Alexander, Homfly and Kauffman polynomials.  We will use this skein
relation to give two new descriptions of the weight systems for Milnor's invariants:  first by means of a
recursive relation which can (like a skein relation) be used to compute the weight system directly from a
chord diagram, and secondly by a characterization in terms of the intersection graph of the chord diagram
\cite{me}.  Along the way, we will give a more direct combinatorial proof that the $\mu$-invariants are finite
type.

We will begin with a brief review of the preliminaries - finite type invariants, Milnor's invariants and
Polyak's skein relation - in Section~\ref{S:prelim}.  In this section we will also give our new proof that
Milnor's invariants are finite type.  In Section~\ref{S:weight} we will give a recursive formula for computing
the Milnor weight systems and we will show that these weight systems detect the presence of a
certain kind of spanning tree in the intersection graph of the chord diagram (this characterization is similar
to Polyak's Gauss formula for the $\mu$-invariants \cite{po}).

\section{Preliminaries} \label{S:prelim}

\subsection{String Links and Link Homotopy} \label{SS:stringlink}

A {\it string link} of $k$ components is a proper embedding of $k$ disjoint line segments into a
solid cylinder so that each strand has one endpoint on the bottom and one on the top, and the order of the
strands is the same at both ends.  More formally, we have the following definition:

\begin{defn} (Habegger and Lin \cite{hl}) \label{D:string}
Let D be the unit disk in the plane and let I = [0,1] be the unit interval.  Choose k points $p_1,..., p_k$ in
the interior of D, aligned in order along the the x-axis.  A {\bf string link} L of k components is a
smooth proper imbedding of k disjoint copies of I into $D \times I$:
$$L:\ \bigsqcup_{i=1}^k{I_i} \rightarrow D \times I$$
such that $L|_{I_i}(0) = p_i \times 0$ and $L|_{I_i}(1) = p_i \times 1$. The image of $I_i$ is
called the ith string of the string link L.
\end{defn}

Note that any string link can be closed up to a link by joining the top and bottom of each component by an arc
which lies outside of $D \times I$ (and which is unlinked with the other arcs).  An isotopy invariant of
string links is a map from the space of string links to some set which is invariant under ambient isotopies of
the embedding of the line segments in the cylinder which leave the endpoints of the line segments fixed.  A
{\it homotopy} invariant is invariant under ambient isotopies, and also under crossing changes of a strand
with itself (but {\it not} under crossing changes between two different strands).

We will also consider the more general notion of a $(l, 2k-l)-tangle$, which has $k$ strands embedded in the cylinder $D \times I$ with $l$ endpoints arranged along the bottom of the cylinder and $2k-l$ along the top.  So the order of the endpoints may be different on the top and bottom, or there may be strands which have both endpoints at the same end of the cylinder.

\subsection{Milnor Invariants} \label{SS:milnor}

Roughly speaking, Milnor's invariants detect how deeply the longitudes of the components of a string link lie
in the lower central series of the link group.

Given a string link $L$, its link group $\pi_1((D^2 \times I) - L)$ has a Wirtinger presentation, generated by
the arcs of the string link diagram.  We also have a presentation of the link group modulo the $q$th subgroup
in its lower central series (see \cite{hl,mi}):

$$\pi_1((D^2 \times I)-L)/(\pi_1)_q((D^2 \times I)-L) = <m_i\ |\ m_il_im_i^{-1}l_i^{-1} = 1, A_q>$$

\noindent where the generators are the meridians $m_i$ of the components of the link, the $l_i$ denote
the longitudes of the components of the link, and $A_q$ denotes the $qth$ subgroup in the lower
central series of the free group on $\{m_i\}$.  So each longitude (and the generators of the
Wirtinger presentation) can be written in $\pi_1/(\pi_1)_q$ as a word in the $m_i$'s.  We
look at the Magnus expansion of the longitudes, replacing each $m_i$ with $1+K_i$ and each
$m_i^{-1}$ with $1-K_i+K_i^2-...$.  We define $\mu_{i_1...i_n,j}(L)$ as the coefficient of
$K_{i_1}...K_{i_n}$ in the Magnus expansion of $l_j$.  This is a well-defined invariant of string links up to
concordance, as long as $q>n$ (it is otherwise independent of $q$).  If the indices $i_1,...i_n,j$ are all
distinct, it is in fact an invariant of string link homotopy.  This definition can easily be extended to
oriented $(l, 2k-l)$ tangles without closed components \cite{po}.

Note in particular that if component $j$ overcrosses all other components, then $\mu_{i_1...i_n,j}$ is
trivial, since the longitude  $l_j$ will be the empty word.  Also, if any component $i_1,...,i_n$ {\it
undercrosses} the other components, then its meridian will never appear in $l_j$, and $\mu_{i_1...i_n,j}$ will
likewise be trivial.

\begin{figure} [h]
    $$\includegraphics{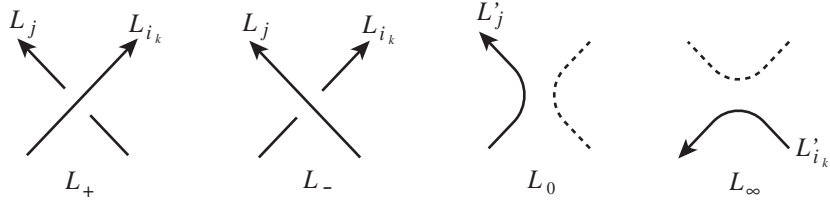}$$
    \caption{Splitting a crossing} \label{F:muskein}
    \end{figure}

However, this definition can be rather unwieldy for computation.  Polyak \cite{po} has shown that
the $\mu$-invariants without repeating indices (i.e. the homotopy invariants) satisfy a skein relation,
similar to the skein relations for the Conway-Alexander and Jones polynomials, which provides another method of
computation.  Consider two string links (or oriented tangles) $L_+$ and $L_-$ which are identical except for
one crossing of components $i_k$ and $j$, where $L_+$ has the positive crossing and $L_-$ has the negative
crossing.  We can "split" this crossing in two ways to define two other oriented tangles.  Let $L'_{i_k}$ be a
strand which follows $L_j$ (respecting orientation) until reaching the crossing, and then follows $L_{i_k}$
against the orientation.  Let $L'_j$ be a strand which follows $L_{i_k}$ and then switches to $L_j$,
respecting the orientation on both components.  Then we define
$L_0 = (L - L_{i_k} - L_j) \cup L'_j$ and $L_\infty = (L - L_{i_k} - L_j) \cup L'_{i_k}$ as shown in
Figure~\ref{F:muskein}.
Polyak showed that Milnor's $\mu$-invariants satisfy the following skein relation:
$$\mu_{i_1...i_k...i_n,j}(L_+) - \mu_{i_1...i_k...i_n,j}(L_-) = \mu_{i_1...i_{k-1},i_k}(L_\infty) \cdot
\mu_{i_{k+1}...i_n,j}(L_0)$$
In particular, in the cases when $k = 1$ or $k = n$, we have:
$$\mu_{i_k...i_n,j}(L_+) - \mu_{i_k...i_n,j}(L_-) = \mu_{i_{k+1}...i_n,j}(L_0)$$
$$\mu_{i_1...i_k,j}(L_+) - \mu_{i_1...i_k,j}(L_-) = \mu_{i_1...i_{k-1},i_k}(L_\infty)$$
$$\mu_{i_k,j}(L_+) - \mu_{i_k,j}(L_-) = 1$$
So we can reduce any computation of $\mu$-invariants to a computation of linking numbers for tangles, which
can be defined (using Gauss's formula) as one-half the sum of the signs of the crossings (a positive crossing
contributes +1, and a negative crossing -1).  Since tangles may have an odd number of crossings, this will not
always yield an integer, as it does for links and string links.

{\sc Notation:}  We will often consider the invariant $\mu_{1...n,n+1}$, or other invariants where the indices
are a sequence of consecutive integers.  For convenience, we will use $\mu_{i;j}$ to denote
$\mu_{i,i+1,...,j-1,j}$.  

\subsection{Finite Type Invariants, Chord Diagrams and Weight Systems} \label{SS:finitetype}

Our treatment of finite type invariants will follow the combinatorial approach of Birman and Lin~\cite{bl}. 
We will give a brief overview of this combinatorial theory; for more details, see \cite{bn,bl}.

\begin{figure} [h]
    $$\includegraphics{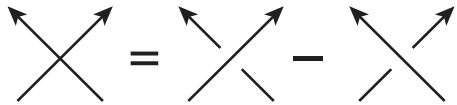}$$
    \caption{Extending invariants to singular links} \label{F:extend}
    \end{figure}

We first note that we can extend any tangle invariant to an invariant of {\it singular} tangles, where a singular $(l, 2k-l)$-tangle is an immersion of a disjoint union of line segments in a solid cylinder (with
$l$ endpoints on the bottom and $2k-l$ on the top) which is an embedding except for a finite number of
isolated double points. We extend a tangle invariant $v$ to singular tangles via the relation in Figure~\ref{F:extend}

An invariant $v$ of singular tangles is then said to be of {\it finite type}, specifically of {\it type
n}, if $v$ is zero on any tangle with more than $n$ double points (where $n$ is a finite nonnegative
integer). We denote by $V_n$ the vector space over ${\mathbb C}$ generated by (framing-independent) finite
type invariants of type $n$.  We can completely understand the space of finite type invariants by
understanding all of the vector spaces $V_n/V_{n-1}$.  An element of this vector space is completely
determined by its behavior on tangles with exactly $n$ singular points.  In addition, since such an
element is zero on tangles with more than $n$ singular points, any other (non-singular) crossing of the
tangle can be changed without affecting the value of the invariant.  This means that elements of
$V_n/V_{n-1}$ can be viewed as functionals on the space of {\it chord diagrams}:

\begin{defn}
A {\bf (l, 2k-l)-tangle chord diagram of degree n} is a disjoint union of k oriented line segments (called the {\bf components} of the diagram) in a cylinder with l endpoints along the bottom of the cylinder and 2k-l along the top, together with $n$ chords (line segments with endpoints on the oriented line segments), such that all of the $2n$ endpoints of the chords are distinct.  The oriented line segments represent a (l, 2k-l)-tangle and the endpoints of a chord represent 2 points identified by the immersion of this tangle into 3-space.  The diagram is determined by the orders of the endpoints on each component.  In the special case of a string link, the components of the diagram are k parallel line segments.
\end{defn}

Functionals on the space of chord diagrams which are derived from finite type tangle invariants will
satisfy certain relations.  This leads us to the definition of a {\it weight system}:

\begin{defn} \label{D:weight}
A {\bf weight system of degree n} is a linear functional $W$ on the space of chord diagrams of degree $n$
(with values in an associative commutative ring ${\bf K}$ with unity) which satisfies the 1-term, 4-term and
antisymmetry relations, shown in Figure~\ref{F:4-term}.  The three arcs of the 4-term relation may belong to
the same or different components.
    \begin{figure} [ht]
    (1-term relation) $$\includegraphics{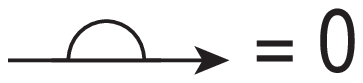}$$
    (4-term relation) $$\includegraphics{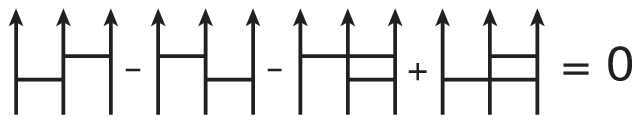}$$
    (antisymmetry relation) $$\includegraphics{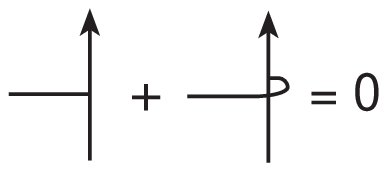}$$
    \caption{The 1-term, 4-term and antisymmetry relations.  No other chords have endpoints on the arcs shown. 
In the 4-term and antisymmetry relations, all other chords of the diagrams are the same.} \label{F:4-term}
    \end{figure}
\end{defn}

The natural map from elements of $V_n/V_{n-1}$ to functionals on chord diagrams is a homomorphism into the
space of weight systems \cite{bn,bl}.  Kontsevich proved the much more difficult fact that
these spaces are isomorphic \cite{bn,ko} (the inverse map is the famous {\it Kontsevich integral}).  So we
can take the dual approach, and simply study the space of chord diagrams of degree $n$ modulo the 1-term and
4-term relations.

When we are looking at {\it homotopy} invariants, such as the $\mu$-invariants, we can add another relation. 
Any homotopy invariant will vanish on a singular tangle with a double point where a component crosses
itself, since the two resolutions differ only by a crossing change of the component with itself.  This means
that the associated weight system will vanish on any chord diagram that has a chord with both endpoints on the
same component.

\subsection{Milnor Invariants are Finite Type} \label{SS:milnorfinite}

It is well-known that Milnor's invariants are finite type invariants of string links \cite{bn2, li}, but the
proofs are rather complicated.  We will use Polyak's skein relation to give a simple combinatorial proof for
the invariants without repeating indices.

\begin{thm} \label{T:mufinite}
The string link homotopy invariant $\mu_{i_1...i_n,j}$ is type $n$.
\end{thm}
{\sc Proof:}  The proof is by induction on $n$.  It is well-known (and easy to prove) that the linking numbers
$\mu_{i,j}$ are type 1.  So we assume inductively that $\mu_{i_1...i_k,j}$ is type $k$ for $k < n$, and
consider a string link $L$ with $n+1$ double points.  Without loss of generality, we can consider the
invariant $\mu_{1...n,n+1}$, denoted $\mu_{1;n+1}$ (this is just a relabeling of the components).  If none of
the double points involves component $n+1$, we can look at the result of bringing component $n+1$ to the front
(so it overcrosses all other components) by a series of crossing changes.  The resulting string link $L'$ has
trivial $\mu_{1...n,n+1}$, and the difference between $L$ and $L'$ is a linear combination of string links
with $n+2$ double points, one of which involves component $n+1$.  So it suffices to show that
$\mu_{1...n,n+1}$ is trivial on any string link with at least $n+1$ double points, at least one of which
involves component $n+1$.  And any such link is a linear combination of string links with exactly $n+1$ double
points (making sure to keep a double point involving component $n+1$).

So we may assume that $L$ has exactly $n+1$ singular crossings $c_1, c_2,...,c_{n+1}$, and that crossing
$c_{n+1}$ involves components $n+1$ and $k$.  Let $L_{e_1,...,e_{n+1}}$, where $e_i = \pm 1$, be the string
link resulting from resolving each singular crossing $c_i$ into a real crossing with sign $e_i$.  Then:

\begin{eqnarray*}
\mu_{1;n+1}(L) & = \sum_{\{e_1,...,e_{n+1}\}} & {(-1)^{e_1...e_{n+1}}\mu_{1;n+1}(L_{e_1,...,e_{n+1}})} \\
 & = \sum_{\{e_1,...,e_n\}} & {(-1)^{e_1...e_n}(\mu_{1;n+1}(L_{e_1,...,e_n,+}) -
\mu_{1;n+1}(L_{e_1,...,e_n,-}))} \\
 & = \sum_{\{e_1,...,e_n\}} & {(-1)^{e_1...e_n}(\mu_{1;k}(L_{e_1,...,e_n,\infty}) \cdot
\mu_{k+1;n+1}(L_{e_1,...,e_n,0}))} \\
 & = \sum_{\{e_1,...,e_{n-1}\}} & {(-1)^{e_1...e_{n-1}}(\mu_{1;k}(L_{e_1,...,+,\infty}) \cdot
          \mu_{k+1;n+1}(L_{e_1,...,+,0})} \\
 & & - \mu_{1;k}(L_{e_1,...,-,\infty}) \cdot \mu_{k+1;n+1}(L_{e_1,...,-,0})) \\
 & = \sum_{\{e_1,...,e_{n-1}\}} & {(-1)^{e_1...e_{n-1}}(\mu_{1;k}(L_{e_1,...,+,\infty}) \cdot
\mu_{k+1;n+1}(L_{e_1,...,+,0})} \\
 & &  - \mu_{1;k}(L_{e_1,...,-,\infty}) \cdot \mu_{k+1;n+1}(L_{e_1,...,+,0}) +
\mu_{1;k}(L_{e_1,...,-,\infty}) \cdot \mu_{k+1;n+1}(L_{e_1,...,+,0}) \\
 & &  - \mu_{1;k}(L_{e_1,...,-,\infty}) \cdot \mu_{k+1;n+1}(L_{e_1,...,-,0})) \\
 & = \sum_{\{e_1,...,e_{n-1}\}} & {(-1)^{e_1...e_{n-1}}(\mu_{1;k}(L_{e_1,...,*,\infty}) \cdot
\mu_{k+1;n+1}(L_{e_1,...,+,0})} \\
 & & + \mu_{1;k}(L_{e_1,...,-,\infty}) \cdot \mu_{k+1;n+1}(L_{e_1,...,*,0}))
\end{eqnarray*}
Where the * indicates that the crossing has become a double point again.  By continuing in this way, we will
eventually obtain:
$$\mu_{1;n+1}(L) = \sum_{r = 0}^n {\sum_{S \subset \{1,...,n\}, |S|=r}{\mu_{1;k}(L_{S,\infty}) \cdot
\mu_{k+1;n+1}(L_{S,0})}}$$
where $L_{S,\infty}$ has $c_i$ singular if $i \in S$ and negative otherwise, and $L_{S,0}$ has $c_i$ positive
if $i \in S$ and singular otherwise.

By our inductive hypothesis, the terms of this sum will be trivial if $r > k-1$ or $n-r > n-k$.  So the only
non-trivial terms are when $r \leq k-1$ and $r \geq k$.  But this is impossible, so {\it all} terms of the sum
will be trivial, concluding the proof.  $\Box$

\section{Milnor Weight Systems} \label{S:weight}

Knowing that the $\mu$-invariants are finite type, it is natural to investigate the associated weight
systems.  We will denote the weight system associated to the invariant $\mu_{i_1...i_n,j}$ by
$M_{i_1...i_n,j}$.  In this section we will use Polyak's skein relation to give a recursive formula for these
weight systems, and use this formula to show that the weight systems detect the presence of a certain kind of
spanning tree in the intersection graph of the chord diagram.  For the remainder of this section, we
will consider the weight system $M_{1;n+1}$ associated with $\mu_{1;n+1}$, since for any chord diagram
$D$ of degree $n$, $M_{i_1...i_n,j}(D) = M_{1;n+1}(D')$, where $D'$ is obtained from $D$ by relabeling the
components.  

\subsection{Recursive formula for $M_{1...n,n+1}$} \label{SS:recursive}

We begin by making some useful observations about the {\it connection graph} of the chord
diagram.

\begin{defn} \label{D:connection}
Given a chord diagram D the {\bf connection graph of D} is the multigraph C(D) whose vertices are the
components of D, and in which vertices i and j are adjacent if and only if there is a chord in D between
components i and j (multiple chords between two components are represented by multiple edges in the connection
graph, chords with both endpoints on the same component are represented by loops).  The number of edges in C(D)
is equal to the degree of D.
\end{defn}

\begin{lem} \label{L:connectiontree}
If D is a chord diagram of degree n on components 1, 2,..., n+1 such that C(D) is {\bf not} a tree, then
$M_{1;n+1}(D) = 0$.
\end{lem}
{\sc Proof:}  $C(D)$ is a graph with $n+1$ vertices and $n$ edges.  If it is not a tree, then the graph must
not be connected (since any connected graph with $n+1$ vertices and $n$ edges is a tree).  So $D$ consists of
two disjoint chord diagrams $D_1$ and $D_2$, where $D_1$ corresponds to the connected component of $C(D)$
containing vertex $n+1$, and $D_2$ is the remainder of the diagram.  Recall that $M_{1;n+1}(D) =
\mu_{1;n+1}(S(D))$, where $S(D)$ is {\it any} singular string link representing the chord diagram $D$ (this is
well-defined since $\mu_{1;n+1}$ is type $n$).  Choose $S(D)$ so that the components of
$S(D_1)$ are everywhere above the components of $S(D_2)$; this is possible since there are no double points
joining these sets of components.  Then every term in the linear combination of string links represented by
$S(D)$ has the components of $D_1$ lying above the components of $D_2$, so the longitude of component $n+1$
will not detect the components in $D_2$, and the word $l_{n+1}$ will not contain any of the meridians of the
components in $D_1$.  Therefore, $\mu_{1;n+1}$ is trivial on every one of these terms, and thus on $S(D)$. 
We conclude that $M_{1;n+1}(D) = 0$, as desired.  $\Box$

In particular, this means that (if $M_{1;n+1}(D)$ is non-trivial) $C(D)$ can have no multiple edges or
loops, so $D$ cannot have multiple chords between two components, or chords with both endpoints on the same
component.  We should also observe that if $D$ has a chord with an endpoint on a component $r > n+1$, then
$M_{1;n+1}(D) = 0$, since $\mu_{1;n+1}$ cannot detect the difference in the two resolutions of the
corresponding double point.

\begin{figure} [h]
    $$\includegraphics{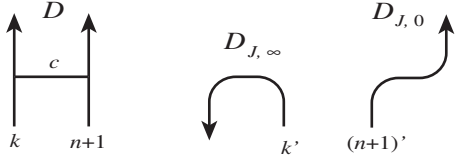}$$
    \caption{Splitting a chord} \label{F:weightskein}
    \end{figure}

Recall from the proof of Theorem \ref{T:mufinite} that
$$\mu_{1;n+1}(L) = \sum_{r = 0}^n {\sum_{S \subset \{1,...,n\}, |S|=r}{\mu_{1;k}(L_{S,\infty}) \cdot
\mu_{k+1;n+1}(L_{S,0})}}$$
On the level of weight systems, this means that
$$M_{1;n+1}(D) = \sum_{r = 0}^{n-1} {\sum_{J \subset D-c, |J|=r}{M_{1;k}(D_{J,\infty}) \cdot
M_{k+1;n+1}(D_{J,0})}}$$
where $c$ is the chord between components $k$ and $n+1$, $J$ is a subset of the chords of $D-c$,
$D_{J,\infty}$ is the chord diagram formed from $D$ by replacing components $k$ and $n+1$ by a new component
$k$ as in Figure~\ref{F:weightskein} and including the chords in $J$, and $D_{J,0}$ is the chord diagram
formed from $D$ by replacing components $k$ and $n+1$ by a new component $n+1$ as in Figure~\ref{F:weightskein}
and including the chords in $(D-c)-J$.

The chords in $J$ must have both endpoints on components in $\{1,...,k,n+1\}$, and the chords in
$(D-c)-J$ must have both endpoints on components in $\{k,k+1,...,n,n+1\}$.  Since $c$ was the only chord
between components $k$ and $n+1$ (by Lemma \ref{L:connectiontree}, there are no multiple chords between pairs
of components), the chords in $J$ must have at least one endpoint on the chords $\{1,...,k-1\}$, and the
chords in $(D-c)-J$ must have at least one endpoint on the chords $\{k+1,...,n\}$.  So, in fact, the choice of
$J$ is determined, and we have:

\begin{thm} \label{T:weightskein}
If D is a chord diagram on n+1 components with a chord c between components k and n+1, then:
\begin{equation} \label{E:skein}
   M_{1;n+1}(D) = M_{1;k}(D_{J,\infty}) \cdot M_{k+1;n+1}(D_{J,0})
\end{equation}
where J is the set of chords in D with at least one endpoint on components 1,..., k-1.  In particular, if k =
1 we have:
$$M_{1;n+1}(D) = M_{2;n+1}(D_{J,0})$$
and if k = n we have:
$$M_{1;n+1}(D) = M_{1;n}(D_{J,\infty})$$
and if $D_{ij}$ is the diagram with a single chord connecting components i and j, then $M_{i,j}(D_{ij}) = 1$.
\end{thm}

\begin{figure} [h]
    $$\includegraphics{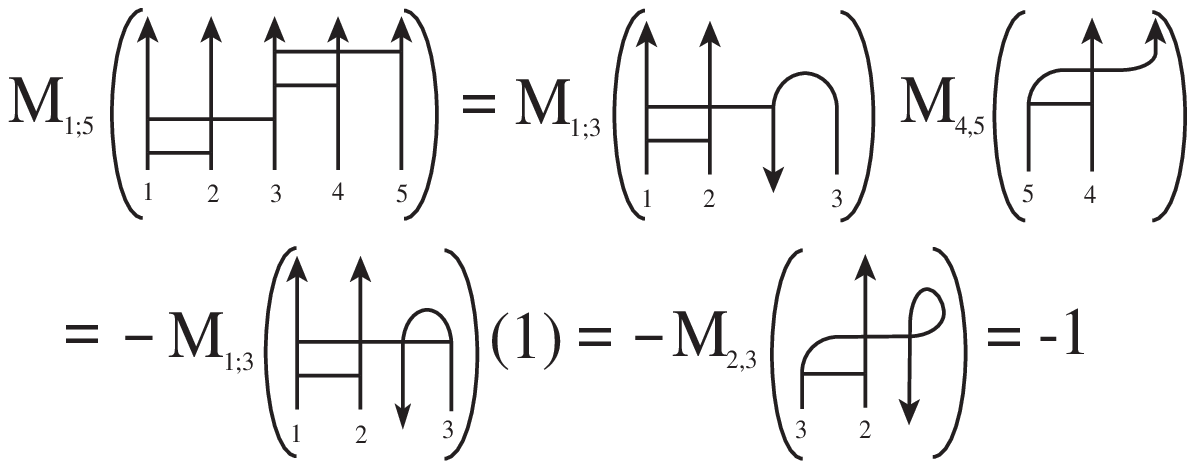}$$
    \caption{Computing the Milnor weight system} \label{F:example}
    \end{figure}

In particular, this product is trivial unless $|J| = k-1$, by Theorem \ref{T:mufinite}.  The product is also
trivial if any chord in $D_{J,\infty}$ or $D_{J,0}$ has an endpoint on one of the "erased" sections of
components $k$ or $n+1$.  Figure~\ref{F:example} shows how to use relation (\ref{E:skein}) to compute
$M_{1;n+1}(D)$.

We will end this section with a useful lemma which uses Theorem \ref{T:weightskein} to give another sufficient
condition for $M_{1;n+1}$ to be trivial.

\begin{lem} \label{L:interlace}
Say that a chord diagram D has chords a and b such that a has endpoints on components i and j, and b has
endpoints on components k and l.  If $i < k < j < l$ (the chords "interlace") then $M_{1;n+1}(D) = 0$.
\end{lem}
{\sc Proof:}  Our proof will be by induction on $n$.  Since $i, j, k, l$ are all distinct, $D$ must have at
least 4 components, so our base case is $n = 3$ ($n+1 = 4$).  In this case $i = 1$, $k = 2$, $j = 3$ and $l =
n+1 = 4$.  Applying relation (\ref{E:skein}) to $D$ using chord $b$, we get $M_{1;4}(D) =
M_{1,2}(D_{J,\infty}) \cdot M_{3,4}(D_{J,0})$, where $J$ is the set of chords with at least one endpoint on
component 1.  So chord $a$ is in $J$.  But since the other endpoint of $a$ is on component 3,
$M_{1,2}(D_{J,\infty}) = 0$, and so $M_{1;4}(D) = 0$.

Inductively, we assume the lemma is true for $n < r$, and consider the case $n = r$.  If $l = n+1 = r+1$, we
can apply relation (\ref{E:skein}) to chord $b$ and find that $M_{1;n+1}(D) = 0$ as in the
base case.  If $l < n+1$, then there is a third chord $c$ with endpoints on components $m$ and $n+1$
(otherwise $C(D)$ is not connected, and $M_{1;n+1}(D) = 0$ as in Lemma \ref{L:connectiontree}).  We apply
relation (\ref{E:skein}) to chord $c$.  If $m > l$ then $D_{J,\infty}$ contains both chords $a$ and $b$, so
$M_{1;m}(D_{J,\infty}) = 0$ by the inductive hypothesis.  Similarly, if $m < i$ then $D_{J,0}$ contains both
$a$ and $b$, so $M_{m+1;n+1}(D_{J,\infty}) = 0$.  Otherwise we must have $i < m < j$ or $k < m < l$, in either
case one of the chords $a$ or $b$ interlaces with chord $c$, so $M_{1;m}(D_{J,\infty}) = 0$ as in the base
case.  This completes the induction. $\Box$

\subsection{Dependence on intersection graphs} \label{SS:IGdependence}

We will show that the weight system $M_{1;n+1}$ depends only on the intersection graph of the chord diagram. 
The intersection graph for a tangle chord diagram is defined as follows:

\begin{defn} \cite{me} \label{D:IGstring}
Let $D$ be a chord diagram with $k$ components (oriented arcs, colored from 1 to $k$) and $n$ chords.  The {\bf intersection graph} $\Gamma(D)$ is the labeled, directed multigraph such that:
\begin{itemize}
    \item $\Gamma(D)$ has a vertex for each chord of $D$.  Each vertex is labeled
by an unordered pair $\{i,j\}$, where $i$ and $j$ are the labels of the components
on which the endpoints of the chord lie.
    \item There is a directed edge from a vertex $v_1$ to a vertex $v_2$ for each
pair $(e_1, e_2)$ where $e_1$ is an endpoint of the chord associated to $v_1$,
$e_2$ is an endpoint of the chord associated to $v_2$, $e_1$ and $e_2$ lie on the
same component of $D$, and the orientation of the component runs from $e_1$ to
$e_2$ (so if the components are all oriented upwards, $e_1$ is below $e_2$).  We
count these edges "mod 2", so if two vertices are connected by two directed edges with the same direction, the edges cancel each other.  If two vertices are connected by a directed edge in each direction, we will simply connect them
by an undirected edge.
\end{itemize}
\end{defn}

\begin{figure} [h]
    $$\includegraphics{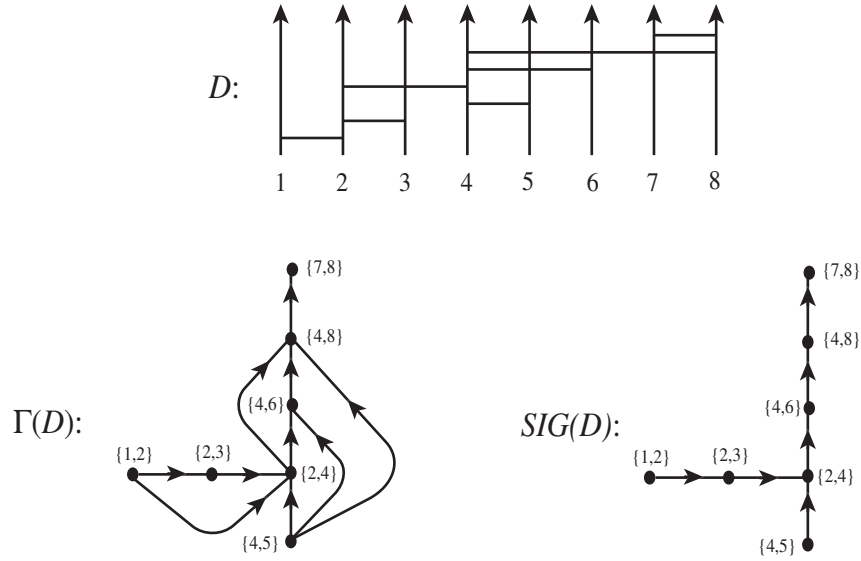}$$
    \caption{The intersection graph and simplified intersection graph of a chord diagram} \label{F:SIG}
    \end{figure}

Figure~\ref{F:SIG} shows an example of a chord diagram and its intersection graph.
Observe that the intersection graph $\Gamma(D)$ determines the connection graph $C(D)$, since $C(D)$ depends
only on the labels of the vertices in $\Gamma(D)$.  We will need the following lemma:

\begin{lem} \label{L:IGCconntree}
If $D_1$ and $D_2$ are chord diagrams on the same tangles such that $\Gamma(D_1) = \Gamma(D_2) = \Gamma$ and $C(D_1) = C(D_2) = C$
is a tree, then $D_1 = D_2$.  (The diagrams are actually identical, not just equivalent modulo the 4-term and
1-term relations.)
\end{lem}
{\sc Proof:}  Since $C$ is a tree, every chord in $D_1$ and $D_2$ has its endpoints on different components,
and two chords can share at most one component.  This means that every edge in $\Gamma$ is directed.  Hence,
the relative positions of the endpoints along each component are determined by the directions of the edges in
$\Gamma$, which completely determines the original chord diagram.  $\Box$

By Lemma \ref{L:connectiontree}, if $\Gamma(D_1) = \Gamma(D_2)$ and $C(D_1) = C(D_2)$ is {\it not} a tree,
then $M_{1;n+1}(D_1) = M_{1;n+1}(D_2) = 0$.  Together with Lemma \ref{L:IGCconntree}, this gives us:

\begin{thm} \label{T:IGdependence}
If $D_1$ and $D_2$ are degree n chord diagrams on the same tangle of n+1 components and $\Gamma(D_1) = \Gamma(D_2)$, then $M_{1;n+1}(D_1) = M_{1;n+1}(D_2)$.
\end{thm}

Using Theorem \ref{T:weightskein}, we could also prove this by showing that we can obtain
$\Gamma(D_{J,\infty})$ and $\Gamma(D_{J,0})$ from $\Gamma(D)$.  This exercise is left to the reader.  This approach also leads to the following useful corollary.

\begin{cor} \label{C:redrawing}
If $D_1$ and $D_2$ are both chord diagrams on tangles with $n+1$ components (but not necessarily the same tangles), and $\Gamma(D_1) = \Gamma(D_2)$, then $M_{1;n+1}(D_1) = M_{1;n+1}(D_2)$.
\end{cor}

\noindent {\sc Proof:}  The result is clearly true for linking numbers - $M_{i,j}(D) = \pm 1$ if and only if there is a single chord connecting components $i$ and $j$, regardless of the tangle.  But since the calculation of $M_{1;n+1}$ can be reduced to computing these linking numbers by Theorem \ref{T:weightskein}, and the reduction can be carried out entirely on the level of intersection graphs, the underlying tangle diagram does not affect the result. $\Box$

As a result of Corollary \ref{C:redrawing}, we can "redraw" chord diagrams as a different diagram with the same intersection graph, without changing the values of the Milnor weight system.  In practice, we will often redraw tangle diagrams as string link diagrams (with disjoint parallel chords) in order to simplify our work.

\subsection{Spanning trees in the intersection graph} \label{SS:IGspantree}

We would like to have a better understanding of how the Milnor weight systems are related to intersection
graphs.  We will show that $M_{1;n+1}$ detects the presence of a certain kind of spanning tree - it is $\pm 1$
if this tree exists, and 0 otherwise.  We will also see how the tree determines the sign.

We will begin by defining a special subgraph of the intersection graph of a chord diagram - the {\it
simplified intersection graph} (or SIG).  The idea of the SIG is to focus on the order of the chords along
each component, and strip away any extraneous edges.

\begin{defn} \label{D:sig}
Let $D$ be a chord diagram with $k$ components and $n$ chords.  The {\bf simplified intersection graph}
$SIG(D)$ is the labeled, directed multigraph such that:
\begin{itemize}
    \item $SIG(D)$ has a vertex for each chord of $D$.  Each vertex is labeled by an unordered pair $\{i,j\}$,
where $i$ and $j$ are the labels of the components on which the endpoints of the chord lie.
    \item There is a directed edge from a vertex $v_1$ to a vertex $v_2$ for each pair $(e_1, e_2)$ where
$e_1$ is an endpoint of the chord associated to $v_1$, $e_2$ is an endpoint of the chord associated to $v_2$,
$e_1$ and $e_2$ lie on the same component of $D$, there are no chords between $e_1$ and $e_2$ on that
component, and the orientation of the component runs from $e_1$ to $e_2$ (so if the components are all
oriented upwards, $e_1$ is {\bf directly} below $e_2$).  We count these edges "mod 2", as we do for the usual
intersection graph.
\end{itemize}
\end{defn}

\noindent $SIG(D)$ is a subgraph of $\Gamma(D)$ which contains all the information needed to reconstruct
$\Gamma(D)$.  Figure \ref{F:SIG} gives an example of a chord diagram with its intersection graph and
simplified intersection graph.

\begin{lem} \label{L:SIGtree}
If $M_{1;n+1}(D) \neq 0$, then SIG(D) is a rooted directed tree.  I.e. SIG(D) is a tree with a unique vertex r
(the root) such that there is a directed path from every other vertex to r.  Moreover, r is labeled \{i, n+1\}
for some i.
\end{lem}
{\sc Proof:}  Let $D$ be a degree $n$ chord diagram on $n+1$ components, with $M_{1;n+1}(D) \neq 0$.  Note that
$SIG(D)$ must be connected - otherwise, either $C(D)$ is disconnected, or there are two chords with endpoints
on the same components (so the edges between the vertices have cancelled "mod 2") and $C(D)$ has a loop. 
Either case contradicts Lemma \ref{L:connectiontree}.

We will first show that $SIG(D)$ is a tree.  Assume $SIG(D)$ contains a loop $v_1v_2...v_kv_1$, where $v_i$ is
labeled by $\{a_i, b_i\}$; so $v_i$ corresponds to a chord $c_{v_i}$ with endpoints on components $a_i$ and
$b_i$.  We first consider the case when all the $a_i$'s are the same - i.e. all the chords have an endpoint on
the same component.  By Lemma \ref{L:connectiontree}, $C(D)$ is a tree, so no two chords
can have the same labels; hence, all the $b_i$'s are distinct.  Then the order of the chords along component
$a = a_i$ is the same as the order of the vertices in the loop $v_1...v_kv_1$, which means that the chord
$c_{v_1}$ is simultaneously above and below the chord $c_{v_k}$ along component $a$.  This is impossible.

But if the $c_{v_i}$'s do not all have endpoints on the same component, then as we follow the loop we move
from one component of $D$ to another, and ultimately return to the original component.  This gives a loop in
$C(D)$, which contradicts Lemma \ref{L:connectiontree}.  So $SIG(D)$ must be a tree.

It remains to show that $SIG(D)$ is a {\it rooted} tree.  We will show that every vertex has outdegree 1
except for one vertex (the root) which has outdegree 0, and that the root must have a label $\{i, n+1\}$. 
Then any directed path from any vertex must ultimately terminate at the root (since there are no loops). 
Assume there is a vertex $v$ (corresponding to chord $c_v$) with outdegree 0 and with label $\{i,j\}$, where
$i < j < n+1$.  If there is a vertex $w$ labeled $\{i,n+1\}$ or $\{j,n+1\}$ then we can apply relation
\ref{E:skein} to the corresponding chord $c_w$.  Since $v$ has outdegree 0, and $SIG(D)$ is connected,
$c_v$ must be above $c_w$ along component $i$ (respectively, $j$).  So after applying relation
\ref{E:skein} to $c_w$, $c_v$ has an endpoint on an erased segment of $i$ (resp. $j$) in $D_{J,0}$ (resp.
$D_{J,\infty}$).  This means that $M_{1;n+1}(D) = 0$, a contradiction.

If we do not have a vertex labeled $\{i,n+1\}$ or $\{j,n+1\}$, there will be some other vertex $u$ labeled
$\{k,n+1\}$, $k \neq i,j$.  The chord $c_u$ does not interlace $c_v$ by Lemma \ref{L:interlace}, so
$c_v$ will be in either $D_{J,0}$ or $D_{J,\infty}$ when we apply relation (\ref{E:skein}) to $c_u$.  By
induction on the number of components, either $M_{1;k}D_{J,\infty} = 0$ or $M_{k+1;n+1}D_{J,0} = 0$; in
either case $M_{1;n+1}(D) = 0$, a contradiction.  So any vertex with outdegree 0 must be labeled $\{i, n+1\}$
for some $i$.  Since there is a directed path through the vertices with label $n+1$, there can be at most one
such vertex; but since $SIG(D)$ is a tree, there must be at least one vertex with outdegree 0.  We conclude
that there is exactly one vertex with outdegree 0, and it has label $n+1$.  This is the root $r$.

We still need to show that no vertex can have an outdegree greater than 1.  Assume vertex $v$ has directed
edges out to vertices $w$ and $u$.  Then we can start at $v$ and follow a directed path through $w$ or through
$u$.  These paths must share some other vertex - if nothing else, they must both terminate at the unique root
$r$.  This would mean that $SIG(D)$ has a (undirected) loop, which is a contradiction.  Therefore, except for
$r$, every vertex has outdegree 1, and $SIG(D)$ is a rooted directed tree.  $\Box$

While it is necessary that $SIG(D)$ be a rooted tree for $M_{1;n+1}(D)$ to be non-trivial, it is not
sufficient.  In the case when $SIG(D)$ is a rooted tree, we can define another spanning tree of $\Gamma(D)$
which we call the {\it branched simplified intersection graph} (BSIG).  Unlike $SIG(D)$, the graph $BSIG(D)$
is dependent on the labels of the vertices (just as the Milnor invariants are dependent on the order of the
indices).  Before we define this new spanning tree, we will look at some implications of the simplified
intersection graph being a tree.

For each component $i$ in $D$, let $r_i$ be the top chord on component $i$ in $D$ (and the corresponding
vertex in $SIG(D)$).  Note that only the root $r$ can be the top chord on two different components.  Since
every other chord has outdegree 1, {\it every} chord is a top chord on one component.  Let $V_i$ be the set of
vertices in $SIG(D)$ with label $\{i,j\}$ for some $j$.  Then we can partition $V_i - \{r_i\}$ into two sets
$L_i$ and $R_i$ (the left and right hand branches), where $L_i = \{\{i,j\}: j < i\}$ and $R_i = \{\{i,j\}: j >
i\}$.  If vertex $\{i,j\}$ is in $L_i$ or $R_i$, then it must be the top chord $r_j$ on component $j$, and
vice versa; so the collection of the $L_i$'s and $R_i$'s is a partition of all the vertices of $SIG(D)$ except
the root $r$.  Obviously, $R_{n+1}$ is empty.

The branched simplified intersection graph takes each path in $SIG(D)$ through the vertices in $V_i$ and
splits it into two paths - one through the vertices in $L_i$ and the other through the vertices in $R_i$ - 
joining the paths at $r_i$.

\begin{defn} \label{D:bsig}
Let $D$ be a chord diagram of degree n with $n+1$ components such that $SIG(D)$ is a directed rooted tree. The
{\bf branched simplified intersection graph} $BSIG(D)$ is the labeled, directed graph such that:
\begin{itemize}
    \item $BSIG(D)$ has a vertex for each chord of $D$.  Each vertex is labeled by an unordered pair $\{i,j\}$,
where $i$ and $j$ are the labels of the components on which the endpoints of the chord lie.
    \item There is a directed edge from a vertex $v$ to a vertex $w$ if both vertices are in $L_i$
(respectively $R_i$) and $w$ is the first vertex in $L_i$ (resp. $R_i$) along the path in SIG(D) from $v$ to
$r_i$.
    \item There are directed edges from the top vertices of $L_i$ and $R_i$ to $r_i$.
\end{itemize}
\end{defn}

\begin{figure} [h]
    $$\includegraphics{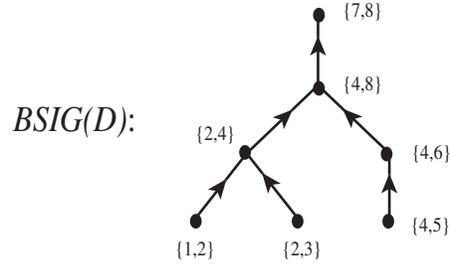}$$
    \caption{The branched simplified intersection graph} \label{F:BSIG}
    \end{figure}

So $BSIG(D)$, when it exists, is also a rooted directed tree which is a spanning tree for $\Gamma(D)$. 
Figure~\ref{F:BSIG} gives the branched simplified intersection graph for the chord diagram in Figure
\ref{F:SIG}.

\begin{defn} \label{D:bsiggood}
We say that BSIG(D) (if it is defined) is {\bf good} if given vertices v with label $\{i,j\}$ and $w$ with
label $\{i,k\}$ we have:
\begin{itemize}
    \item If i = n+1 and $v, w \in L_{n+1}$, then $j < k$ if and only if there is a directed path from v to w.
    \item If $v, w \in L_i$ ($i\neq n+1$), then $j < k$ if and only if there is a directed path from w to v.
    \item If $v, w \in R_i$ ($i\neq n+1$), then $j < k$ if and only if there is a directed path from v to w.
\end{itemize}
\end{defn}

So the graph in Figure \ref{F:BSIG} is a {\it good} branched simplified intersection graph.  We can now state
the main result of this section.

\begin{thm} \label{T:milnorbsig}
If D is a chord diagram on n+1 components such that BSIG(D) exists and is good, then $M_{1;n+1}(D) = (-1)^L$,
where $L = \sum_{i \neq n+1}{|L_i|}$ and $|L_i|$ is the cardinality of $L_i$.  Otherwise, $M_{1;n+1}(D) = 0$.
\end{thm}

For example, if $D$ is the chord diagram in Figure \ref{F:SIG}, with $BSIG(D)$ shown in Figure \ref{F:BSIG},
then $L = |L_2| + |L_4| = 2$, and $M_{1;8}(D) = (-1)^2 = 1$.  We will prove this theorem in two parts.  In
Lemma \ref{L:trivial} we show that if $D$ does not have a good $BSIG(D)$, then $M_{1;n+1}(D) = 0$.  In Lemma
\ref{L:nontrivial} we show that if $D$ does have a good $BSIG(D)$, then $M_{1;n+1}(D) = (-1)^L$.

\begin{lem} \label{L:trivial}
If chord diagram D does not have a good BSIG, then $M_{1;n+1}(D) = 0$.
\end{lem}
{\sc Proof:}  If $BSIG(D)$ does not exist then $SIG(D)$ is not a rooted directed tree, so
$M_{1;n+1}(D) = 0$ by Lemma \ref{L:SIGtree}.  So we can assume that $BSIG(D)$ exists, but that it
fails one of the three conditions in Definition \ref{D:bsiggood}.
\\

\noindent{\it Case 1:}  We first consider the case when $BSIG(D)$ has vertices $v$ and $w$ with labels
$\{j,n+1\}$ and $\{k,n+1\}$, respectively ($j, k < n+1$).  We denote the corresponding chords $c_v$ and
$c_w$.  Since the connection graph is a tree, $j \neq k$.  Assume that $j < k$, but there is a directed
path from $w$ to $v$.  This means that the chord $c_w$ lies below $c_v$ along component $n+1$ in the diagram
$D$.  If we apply relation (\ref{E:skein}) to $c_w$, we have $M_{1;n+1}(D) = M_{1;k}(D_{J,\infty}) \cdot
M_{k+1;n+1}(D_{J,0})$, where $J$ is the set of chords in $D$ with at least one endpoint on components $1,...,
k-1$.  So $c_v$ is in $J$.  But since the other endpoint of $c_v$ is on the erased part of component
$n+1$ in $D_{J,\infty}$, $M_{1;k}(D_{J,\infty}) = 0$, and so $M_{1;n+1}(D) = 0$.  By exactly the same
argument, if $k < j$ and there is a directed path from $v$ to $w$, then $M_{1;n+1}(D) = 0$.  (There must be
some directed path between $v$ and $w$, since one of the corresponding chords is below the other along
component $n+1$.)
\\

\noindent{\it Case 2:}  We now consider the case when $BSIG(D)$ does not satisfy the second condition in
Definition \ref{D:bsiggood}.  Say we have vertices $v$ and $w$ with labels $\{j,i\}$ and $\{k,i\}$,
respectively, and that $v, w \in L_i$ (so $j,k < i < n+1$).  Assume that $j < k$, but there is a directed path
from $v$ to $w$.  Since $D$ is connected, there is some vertex $u$ in $BSIG(D)$ with label $\{m,n+1\}$.  By
Lemma \ref{L:interlace}, if $j < m < i < n+1$ then $M_{1;n+1}(D) = 0$ and we're done.  So assume that either
$m \leq j$ or $m \geq i$.

We first consider when $m = i$.  Applying relation (\ref{E:skein}) to chord $c_u$, we have $M_{1;n+1}(D) =
M_{1;i}(D_{J,\infty}) \cdot M_{i+1;n+1}(D_{J,0})$, where $J$ is the set of chords in $D$ with at least one
endpoint on components $1,..., i-1$.  So both $c_v$ and $c_w$ are in $J$.  We may assume that $c_v$ and $c_w$
are each below $c_u$ along component $i$, or it will have an endpoint on an erased arc and
$M_{1;i}(D_{J,\infty})$ will be trivial.  But then their order along the new component $i$ in $D_{J,\infty}$
will be reversed (see Figure \ref{F:weightskein}), so in $BSIG(D_{J,\infty})$ there is a directed path from
$w$ to $v$.  By Case 1, this means $M_{1;i}(D_{J,\infty}) = 0$, and so $M_{1;n+1}(D) = 0$.

On the other hand, if $m = j$ then applying relation (\ref{E:skein}) to $c_u$ yields $M_{1;n+1}(D) =
M_{1;j}(D_{J,\infty}) \cdot M_{j+1;n+1}(D_{J,0})$, where $J$ is the set of chords in $D$ with at least one
endpoint on components $1,..., j-1$.  In this case, both chords $c_v$ and $c_w$ are in $D_{J,0}$, with $c_v$
now connecting component $i$ and the new component $n+1$.  So the unique directed path from $v$ to the
root of $SIG(D_{J,0})$ passes through vertices with label $n+1$.  But this path must contain $w$, so $v$ and
$w$ have the same label and $C(D)$ is not a tree, which is impossible.

Finally, if $m < j$ or $m > i$ we can proceed inductively on the number of components in the diagram, since
$v$ and $w$ will both be in either $D_{J,\infty}$ or $D_{J,0}$.  So we can apply relation (\ref{E:skein})
until we can use one of the arguments above.  As in Case 1, the other direction of the if and only if is
proved in the same way.
\\

\noindent{\it Case 3:}  Our last case is when $BSIG(D)$ does not satisfy the third condition in Definition
\ref{D:bsiggood}.  Say we have vertices $v$ and $w$ with labels $\{i,j\}$ and $\{i,k\}$, respectively, and
that $v, w \in R_i$ (so $i < j,k < n+1$).  Assume that $j < k$, but there is a directed path from $w$ to $v$. 
Again there is some vertex $u$ with label $\{m,n+1\}$, and we may assume that $m \leq i$ or $m \geq k$.  As in
Case 2, we can inductively reduce the problem to when $m = i$ or $m = k$.  If $m = i$, both $v$ and $w$ will
have endpoints on the new component $n+1$ in $D_{J,0}$ after applying relation (\ref{E:skein}) to chord $c_u$. 
But we will still have $j < k$ and a directed path from $w$ to $v$, which means that $M_{i+1;n+1}(D_{J,0}) =
0$ by the argument in Case 1.

On the other hand, if $m = k$, then $w$ will have a label $n+1$ in $D_{J,\infty}$ after applying relation
(\ref{E:skein}) to $c_u$, and so should be the root along component $i$.  But this is impossible, since there
is a directed path from $w$ to $v$, so we conclude that $M_{1;k}(D_{J,\infty}) = 0$.

We conclude in every case that if $BSIG(D)$ is not good, then $M_{1;n+1}(D) = 0$. $\Box$

\begin{lem} \label{L:nontrivial}
If chord diagram D has a good BSIG, then $M_{1;n+1}(D) = (-1)^L$, where $L = \sum_{i \neq n+1}{|L_i|}$ and
$|L_i|$ is the cardinality of $L_i$.
\end{lem}
{\sc Proof:}  Our proof is by induction on $n$.  If $n = 1$ (so $n+1 = 2$), the only non-trivial chord diagram
with a good $BSIG$ consists of a single chord between components 1 and 2.  Call this diagram $D_{12}$.  Then
$M_{1,2}(D_{12}) = 1$.  Since $L = |L_2| = 0$, $(-1)^L = 1 = M_{1,2}(D_{12})$, so the Lemma is true in this
case.

\begin{figure} [h]
    $$\includegraphics{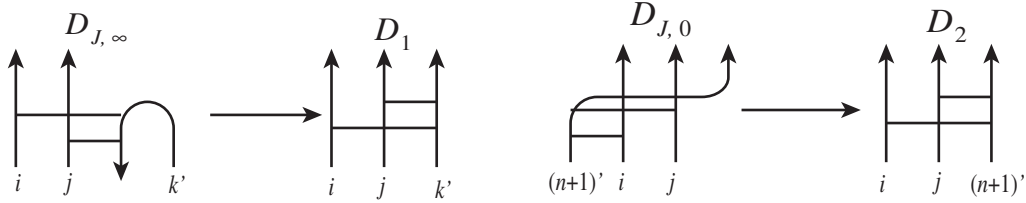}$$
    \caption{Redrawing the chord diagrams} \label{F:straighten}
    \end{figure}

For our inductive step, consider a vertex $v$ in $D$ with label $\{i,n+1\}$ (with
corresponding chord $c_v$).  Choose $v$ to minimize the label $i$ among all vertices with one label $n+1$. 
Applying relation (\ref{E:skein}) to $c_v$, we get $M_{1;n+1}(D) = M_{1;i}(D_{J,\infty}) \cdot
M_{i+1;n+1}(D_{J,0})$, where $J$ is the set of chords in $D$ with at least one endpoint on components $1,...,
i-1$.  So $L_i \subset D_{J,\infty}$ and $R_i \subset D_{J,0}$.  Since we minimized $i$, and $BSIG(D)$ is
good, $c_v$ is the lowest chord on component $n+1$ in $D$.  So all other chords with an endpoint on component
$n+1$ are in $D_{J,0}$.  If we "redraw" the chord diagrams $D_{J,\infty}$ and $D_{J,0}$ in standard form (so
all the components are parallel line segments, arranged in increasing order), as allowed by Corollary \ref{C:redrawing} we obtain diagrams $D_1$ and $D_2$ shown in Figure~\ref{F:straighten}.
By the antisymmetry relation of Definition \ref{D:weight}, $D_{J,\infty} = (-1)^{|L_i|}D_1$ and $D_{J,0} =
D_2$.  We need to show that $D_1$ and $D_2$ have good $BSIG$'s.  Since the order of the chords in $L_i$ has
been reversed in $D_1$, $BSIG(D_1)$ is still good on components $1,...,i$.  

$L_{n+1}(D_2) = L_{n+1}(D) \cup R_i(D)$ (since $c_v$ was the top chord on component $i$), with all the
chords from $R_i(D)$ lying below the chords from $L_{n+1}(D)$.  We just need to show that if there is a vertex
in $L_{n+1}(S)$ with label $\{j,n+1\}$, and a vertex in $R_i(D)$ with label $\{i,k\}$, then $k < j$.  But if
not, the two chords interlace, so $BSIG(D)$ was not good, which is a contradiction. So $BSIG(D_2)$ is good on
components $i+1,...,n+1$.

By the inductive hypothesis, $M_{1;i}(D_1) = (-1)^{\sum_{j < i}{|L_j(D_1)|}} = (-1)^{\sum_{j < i}{|L_j(D)|}}$
and $M_{i+1;n+1}(D_2) = (-1)^{\sum_{i < j < n+1}{|L_j(D_2)|}} = (-1)^{\sum_{i < j < n+1}{|L_j(D)|}}$.  So
$M_{1;n+1}(D) = (-1)^{|L_i|}(-1)^{\sum_{j < i}{|L_j|}}(-1)^{\sum_{i < j < n+1}{|L_j(D)|}} = (-1)^L$.  This
completes the induction and the proof. $\Box$

Combining these lemmas gives the proof of Theorem \ref{T:milnorbsig}.

\small

\normalsize

\end{document}